\documentclass{IEEEtran}

\usepackage{cite}
\usepackage{amsmath,amssymb,amsfonts}
\usepackage{algorithmic}

\usepackage{textcomp}
\def\BibTeX{{\rm B\kern-.05em{\sc i\kern-.025em b}\kern-.08em
    T\kern-.1667em\lower.7ex\hbox{E}\kern-.125emX}}

\IEEEoverridecommandlockouts%% The preceding line is only needed to identify funding in the first footnote. If that is unneeded, please comment it out.
\usepackage{graphicx}
\usepackage{cite}
\usepackage{picinpar}
\usepackage{amsmath}
\usepackage{url}
\usepackage{flushend}
\usepackage[latin1]{inputenc}
\usepackage{colortbl}
\usepackage{soul}
\usepackage{multirow}
\usepackage{pifont}
\usepackage{color}
\usepackage{alltt}
\usepackage[hidelinks]{hyperref}
\usepackage{enumerate}
\usepackage{siunitx}

\usepackage{epstopdf}
\usepackage{pbox}
\usepackage{verbatim}
\usepackage{leftidx}
\usepackage{amssymb}
\usepackage{xcolor}
\usepackage{listings}
\usepackage{textcomp}

\newtheorem{Theorem}{Theorem}[section]
\newtheorem{Proposition}{Proposition}[section]

\newtheorem{Corollary}{Corollary}[section]
\newtheorem{Lemma}[Theorem]{Lemma}
\newtheorem{Definition}{Definition}[section]

 \newcommand{\R}{\mathbb{R}}
  
  \newcommand{\dn}{\mathbf{d}}
  \newcommand{\zero}{0}

\newenvironment{psmallmatrix}
  {\left[\begin{smallmatrix}}
  {\end{smallmatrix}\right]}

\newcommand{\RNum}[1]{\uppercase\expandafter{\romannumeral #1\relax}}
\graphicspath{ {Pictures/} }

\DeclareMathOperator{\rank}{rank}
\def\BibTeX{{\rm B\kern-.05em{\sc i\kern-.025em b}\kern-.08em
    T\kern-.1667em\lower.7ex\hbox{E}\kern-.125emX}}
    
\usepackage{comment}

\begin{document}
\title{	Homogeneous control design using  invariant ellipsoid method} 
\author{Siyuan WANG,  Andrey POLYAKOV, Gang ZHENG, Xubin PING, Driss BOUTAT

%\thanks{This paragraph of the first footnote will contain the date on 
%which you submitted your brief for review. It will also contain support 
%information, including sponsor and financial support acknowledgment. For 
%example, ``This work was supported in part by the U.S. Department of 
%Commerce under Grant BS123456.'' }
\thanks{Siyuan WANG and Driss Boutat  are with INSA Centre Val de Loire, Bourges, France (e-mail: Siyuan.wang@insa-cvl.fr, Driss.Boutat@insa-cvl.fr). }
\thanks{Andrey POLYAKOV and Gang ZHENG are with INRIA Lille, Lille, France. (e-mail: andrey.polyakov@inria.fr, gang.zheng@inria.fr).}
\thanks{Xubin PING is with 
Xidian University, Xian, China (e-mail: pingxubin@126.com).}
}

%\author{\IEEEauthorblockN{1\textsuperscript{st} Siyuan WANG}
%\IEEEauthorblockA{\textit{PRISME } \\
%\textit{ INSA Centre Val de Loire}\\
%Bourges, France \\
%siyuan.wang@insa-cvl.fr}
%\and
%\IEEEauthorblockN{2\textsuperscript{nd} Andrey POLYAKOV}
%\IEEEauthorblockA{\textit{Valse} \\
%\textit{Centre Inria - Lille Nord Europe}\\
%Lille,France \\
%andrey.polyakov@inria.fr}
%\and
%\IEEEauthorblockN{3\textsuperscript{rd} Gang ZHENG}
%\IEEEauthorblockA{\textit{ Defrost} \\
%\textit{Centre Inria - Lille Nord Europe}\\
%Lille, France \\
%gang.zheng@inria.fr}
%\and
%\IEEEauthorblockN{4\textsuperscript{rd} Driss BOUTAT}
%\IEEEauthorblockA{\textit{ PRISME} \\
%\textit{INSA Centre Val de Loire}\\
%Bourges, France \\
%driss.boutat@insa-cvl.fr}
%\and
%\IEEEauthorblockN{5\textsuperscript{rd} Xubin Ping}
%\IEEEauthorblockA{\textit{ ?????} \\
%\textit{Xidian University}\\
%Xian, China 5\\
%pingxubin@126.com}
%}

\maketitle

\begin{abstract}
The invariant ellipsoid method is aimed at minimization of the smallest invariant and attractive set of a linear control system operating under bounded external  disturbances.  This paper extends this technique to a class of  the so-called generalized  homogeneous system by defining the $\dn-$homogeneous invariant/attractive ellipsoid. The generalized homogeneous optimal (in the sense of invariant ellipsoid) controller allows further improvement of the control system providing a faster convergence and better precision. Theoretical results are supported by numerical simulations and experiments.
%The problem of homogeneous controller design based on the invariant set method  is addressed. The equivalent condition of $\dn$-homogeneous ball being invariant set is proposed.  The method of invariant set reduces the design of optimal  controller to the problem of looking for the  smallest invariant set of the closed loop system. Based on the optimized linear controller and minimal invariant set, an optimal dilation is applied to the    homogeneous controller which guarantees that the $\dn$-homogeneous  invariant set is the same to the linear case. This method is easy to implement and meaningful for the practice. The  performance of proposed approach is supported by the simulation results.
%The problem of designing a globally bounded homogeneous controller  by the optimized linear controller is addressed. The method of invariant set reduces the design of optimal  controller to find the  smallest invariant set of the closed loop system. 
% Based on the optimized linear controller and minimal invariant set, an optimal dilation is applied to design the globally bounded   homogeneous controller which guarantees that the $\dn$-homogeneous  invariant set is the same to linear case. This method is easy to implement and meaningful for the practice. The  performance of proposed approach is supported by the simulation results.
\end{abstract}

\begin{IEEEkeywords}
Homogeneity,  Invariant set, LMI
\end{IEEEkeywords}

\markboth{IEEE Transactions }%
{}

\definecolor{limegreen}{rgb}{0.2, 0.8, 0.2}
\definecolor{forestgreen}{rgb}{0.13, 0.55, 0.13}
\definecolor{greenhtml}{rgb}{0.0, 0.5, 0.0}

%%%%%%%%%%%%%%%%%%%%%%%%%%%%%%%%%%%%%%%%%%%%%%%%%%%%%%%%%%%%%%%%%%%%%%%%%%%%%%%%
\section{Introduction} \label{sec:intro}

During the formulation of any control issue, there is always a discrepancy between the actual system and the mathematical model used for control design. This mismatch comes from the unmodelled dynamics, uncertainties in system parameters or the approximation of complex plant. However, the engineer needs to guarantee that the designed controller is able to achieve the required performance despite of all these mismatches. This leads to  the development of the so-called robust control methods solving this problem.  
Robust control design is an approach dealing with perturbations of a nominal system. Its objective is to attain the certain level of performance or stability in the system despite the presence of bounded disturbance.  Several well-known  methodologies of robust control design have been invented such as sliding mode control, $H_\infty$ approach, attractive/invariant ellipsoid method, and so on.
The sliding mode control (SMC) as a 
robust control design methodology is known since 1960s in Russia. The first survey paper by V. Utkin is published in English in 1977  \cite{utkin1977variable}. Sliding mode methodology has been devised for both linear and nonlinear system \cite{utkin2013sliding}, \cite{shtessel2014sliding} and delivers a good performance  in numerous real-world scenarios \cite{utkin2017sliding}. The concepts of the $H_\infty$ control theory  to solve the robust stabilization problem for linear system can be found in \cite{zames1981feedback,orlov2014advanced,kimura1984robust}. Later the $H_\infty$ control methodology  was extended
to the various systems \cite{gershon2005h}.
%in \cite{francis1987course,khargonekar1988h,khargonekar1990robust} for the linear system with different uncertainty. 

%{\color{green}  
%The  suppression of bounded external disturbances is one of the important problem in robust control theory. For example,  
%linear quadratic gaussian (LQG) problem formulated by R.Kalman deals with  the random gaussian disturbance disturbance. In $H_{\infty}$-optimization, the external noise is regarded as random  or $L_2$-bounded which is decreasing with time. In the practical cases, the external disturbance is usually known as bounded without any other information or property. The problem of suppressing external bounded disturbance has  a long history,  for the discrete system, it is firstly formulated by Yakubovich \cite{yakubovich1975solution} and solved in \cite{vidyasagar1986optimal} for some  special cases. The full solution are given by Barabanov and Dahleh in \cite{barabanov1984optimal,dahleh19871} where this problem  is called $l_1$-optimization. However, one serious drawback of  the optimal  $l_1$ controller is that it can be arbitrarily high order.   This kind of problem is solved by the method of dynamical programming in  \cite{elia2000minimization,glover1971control}. 

%For the continuous time case, suppression of bounded external disturbance is usually considered as a difficult problem in control theory \cite{polyak2005hard}. 

The basic ideas of the attractive/invariant method were introduced in the papers \cite{bertsekas1971recursive,glover1971control,schweppe1968recursive}.
One of the main features of this method is that the states of the robustly stabilized system converge to a minimal (in some sense) ellipsoidal set regardless of perturbations or uncertainties satisfying certain bounds. 
The attractive/invariant ellipsoid method is widely applied  to various control and estimation problems for both linear \cite{khlebnikov2011optimization} and nonlinear plants \cite{Poznyak_etal2014:Book}. The key feature of the method is  by using  Linear Matrix Inequalities (LMIs) for control parameters tuning \cite{boyd1994linear}. For instance,  the problem of suppressing bounded additive external disturbance in terms of the invariant ellipsoid is studied in \cite{khlebnikov2011optimization} for linear plant.  It turns out to be equivalent to a Semi-Definite Programming (SDP) problem minimizing the invariant ellipsoid of the closed-loop linear system to design the optimal robust feedback. The further suppression of bounded external disturbance beyond a certain level   using the same \textit{linear feedback} strategy seems infeasible. The purpose of this paper is to extend the attractive/invariant ellipsoid technique to a class of controller known as  generalized homogeneous controllers \cite{Polyakov2020:Book}, which may provide a better quality of control   such as faster convergence, better robustness and smaller overshoots.

A symmetry with respect to a dilation is referred to as homogeneity \cite{Zubov1958:IVM}, \cite{kawski1990homogeneous}. This concept is applied  in the field of control theory and experiment for the purpose of system analysis, controller and  observer design (e.g. \cite{hermes1986nilpotent,BhatBernstein2005:MCSS,levant2005homogeneity,orlov2004finite,polyakov2016robust,wang2021generalized,wang2021generalizedObserver} and references therein). 
%The reason is that the local stability of homogeneous system leads to a global one, the convergence rate can be estimated by its homogeneous degree and the homogeneous system is robust under external disturbance and time delay.  
The standard (Euler) homogeneity means that  a function $f(x)$ remains invariant with respect to the scaling of its argument $f(e^s x) = e^{\nu s}f(x),\forall s,x\in\R$, where the constant $\nu$ is called the homogeneity degree. The generalized (weighted) dilation of vector $x=(x_1,x_2,..,x_n)^\top \in \R^n$ is introduced by V.I. Zubov in 1958 \cite{Zubov1958:IVM}: $(x_1,x_2,..,x_n) \to (e^{r_1 s}x_1,e^{r_2 s}x_2,...,e^{r_n s}x_n),s\in\R$, where the positive numbers $r_1,r_2,...,r_n$ are the weights specifying the dilation rate of each coordinate. Nonlinear (geometric) dilutions are studied in \cite{kawski1991families}, \cite{Rosier1993:PhD}, \cite{bernuau2013robustness}. This paper deals with the \textit{linear geometric dilation} \cite{polyakov2019sliding} given by $x\to e^{G_\dn s}x$, where $G_\dn \in \R^{n\times n}$ is anti-Hurwitz matrix\footnote{The matrix $G_\dn\in\R^n$ is anti-Hurwitz, if $-G_\dn$ is Hurwitz.}. 
The homogeneity as a relaxation of linearity can provide an extra degree of freedom for further minimization of the disturbance effects \cite{mera2016finite}. 

The key contribution of this article is the extension of the attractive/invariant ellipsoid approach to a category of generalized homogeneous control systems\cite{zimenko2020robust}.
Specifically, we present the concept of a homogeneous invariant ellipsoid through the use of the canonical homogeneous norm \cite{polyakov2019sliding}
and derive its characterization in terms of LMIs. We also show that an optimal homogeneous  control design is obtained by solving an SDP problem minimizing the homogeneous attractive/invariant ellipsoid of the system. As an example, the controlled rotary inverted pendulum is studied, where it demonstrate that the optimal (in the sens of attractive ellipsoid) homogeneous controller can stabilize the pendulum with a  better precision than the optimal (in the same sense) linear controller without any degradation of the control signal quality.

The structure of paper is  as follows:  The problem statement is outlined  in Section \ref{Pro_state}. Section \ref{Preliminary} delves into the concept of  homogeneous theories.  The main results about the conditions of being a $\dn$-homogeneous attractive/invariant ellipsoid and  the optimal homogeneous   controller design  are explained in Section \ref{Pre_results} and Section \ref{Main_results}. The last Section \ref{sim_results} provides the simulation results to support the proposed theories.

{\itshape Notation :}
$\mathbb{R}$:  set of real numbers, 
$\mathbb{R}_+ = \{x\in \mathbb{R}:x> 0\}$;
$\|x\| =\sqrt{x^\top x}  $: norm in $\R^n$;  
 $\|x\|_P = \sqrt{x^\top P x}$:  weight norm in $\mathbb{R}^n$;
$\mathrm{diag}\{\lambda_i\}_{i=1}^n$ : the diagonal matrix with  elements $\lambda_i$;
$P\succ0 (\prec0,\succeq 0,\preceq 0)$ for $P\in \mathbb{R}^{n\times n}$ means that the matrix $P$ is symmetric and positive (negative) definite (semi definite);  
$\lambda_{min}(P)$ and $\lambda_{max}(P)$ : the minimal and maximal eigenvalues  of  matrix $P=P^{\top}$;
 for $P\succeq 0$ the square root of $P$ is a matrix $M=P^{\frac{1}{2}}$ such that  $M^2=P$; $\mathrm{tr}(P)$ :  trace of matrix $P$; if $V:\R^n\mapsto [0,+\infty)$ is a positive definite function and $\dot x=f(t,x)$ is an ODE with $f$ continuous on $c$ then we denote the time derivative of $V$ along all solutions of the ODE as 
 $
 \frac{d}{dt} V(x):=\frac{\partial V}{\partial x}f(t,x);
 $
 a function $\gamma : \R_{\geq 0} \to \R_{\geq 0}$ is said to be of class $\mathcal{K}$ if it is continuous strictly increasing and satisfies $\gamma(0) = 0$; a function $\beta : \R_{\geq 0}\times \R_{\geq 0} \to \R_{\geq 0}$ is said to be of class $\mathcal{KL}$ if for each fixed $t$ the mapping $\beta(\cdot,t)$ is of class $\mathcal{K} $ and for each fixed $s$ it is decreasing to zero on t as $t\to \infty$;  $L^{\infty}$ is the space of Lebesgue measurable essentially bounded function $\sigma : \R_+ \to \R^n$ with norm defined as $\|\sigma\|_{L^\infty} : =\mathrm{ess} \sup_{t\in\R_+}\|\sigma(t)\|_\infty<+\infty$; $1_n\in \R^n$ is a vector with all element is 1 
 % }

\section{Problem statement}\label{Pro_state}

%%%%%%%%%%%%%%%%%%%%%%%%%%%%%%%%%%%%%%%%%%%%%%%%%%%%
%%%%%%%%%%%%%%%%%%%%%%%%%%%%%%%%%%%%%%%%%%%%%%%%%%%

%%%%%%%%%%%%%%%%%%%%%%%%%%%%%%%%%%
%%%%%%%%%%%%%%%%%%%%%%%%%%%%%%%%%%%%
 \label{sec:Pro_statement}
% \subsection{Problem statement}
In this paper, we  deal with the linear time-invariant system with additive perturbations,
\begin{equation}\label{eq:control_sys}
\dot{x}(t) = A x(t) + B u(t) +D\omega(t), \quad t>0
\end{equation}
where $x(t)\in \R^n$ is the system state,  $\omega(t)\in \R^p$ is an external bounded disturbance, $u(t) \in \R^m$ is the control input, $A\in \R^{n\times n}$ is nilpotent, $B\in \R^{n\times m}$,  $D\in \R^{n\times p}$ are constant matrix. We assume that the pair $\{A,B\}$ is controllable.

We study the feedback stabilization problem of the system \eqref{eq:control_sys} with perturbations. In the general case, it is not possible to achieve a precise stabilization of this system to a zero state, hence a feedback controller $u(x)$ needs to be optimized in order to minimize (in some sense) the effect of external disturbance $\omega$. For the case of a linear stabilizing feedback there are several well-known methodologies tackling this problem. For example,  $H_{\infty}$ and $H_2$ algorithms (see e.g., \cite{kwakernaak2002h2,khargonekar1990robust}) suggest to optimize a norm of a transfer function, while the attractive/invariant  ellipsoid method \cite{khlebnikov2011optimization,Poznyak_etal2014:Book} optimizes (in a certain sense) the attractive/invariant set of the system \eqref{eq:control_sys} with a linear feedback.  This paper is aimed at  an extension of the latter approach to a special  class of generalized homogeneous \cite{Polyakov_etal2018:TAC,polyakov2019sliding,zimenko2020robust,wang2021generalized} feedback laws. The mentioned nonlinear control algorithms are known to be efficient for finite/fixed-time stabilization of linear plants. Our objective is to design a generalized  homogeneous controller that minimizes the attractive/invariant ellipsoid of the closed-loop system \eqref{eq:control_sys} with 
\begin{equation}\label{eq:disturbance_assumption}
    \omega^{\top}\!(t) Q\omega(t)\leq 1, \quad \omega\in  L^{\infty}(\R_+,\R^p), 
\end{equation}
where $0\prec Q=Q^{\top}\in \R^{p\times p}$ is a known matrix.
%A linear vector field is generalized homogeneous of non-zero degree if the matrix $A$ is nilpotent \cite{Polyakov2020:Book}.

%We propose a procedure for upgrading a linear controller (optimal in the sense of minimal invariant ellipsoid \cite{khlebnikov2011optimization}) to a generalized homogeneous one preserving the minimal invariant  set. We also investigates possible advantages of the obtained feedback law. 
% This article studies how to further minimize the effect of external bounded disturbance for system \eqref{eq:control_sys} by using the optimized parameters of linear controller and the generalized homogeneous controller. In the following section, we firstly provide the equivalent condition of $\dn$-homogeneous ball $\epsilon_x(\rho,P)$ to be a state invariant set for system \eqref{eq:Auto_system_disturbance}. Next,  the method of designing an homogeneous  controller to have an invariant $\dn$-homogeneous ball $\epsilon_x(\rho,P)$ is introduced. Finally an optimal homogeneous controller is designed by applying the optimized linear controller. %This further improves the precision and stability of system.

\section{Preliminaries: Elements of the Homogeneity Theory}\label{Preliminary}

{
\subsection{Linear dilation and monotonicity}
The homogeneity is a symmetry of an object (e.g. a function or a set) with respect to a group of transformations called dilation. In this paper we deal only with the so-called linear dilation\cite{Polyakov2020:Book}: 
$x\mapsto \dn(s)x, x\in \R^n, s\in \R$,
where $s\in \R$ is a parameter of the dilation and 
$$
\dn(s)=e^{sG_{\dn}}=\sum_{i=0}^{+\infty} \frac{s^iG_{\dn}^i}{i!}.
$$
The anti-Hurwitz matrix $G_{\dn}\in \R^{n\times n}$ is known as \textit{the generator of the dilation $\dn$}.

The monotonicity of dilation holds significance in the characterization of homogeneous geometric structures in $\R^n$ and the analysis of homogeneous control systems..
\begin{Definition}\itshape \cite{Polyakov2020:Book}\label{def:monotone}
Dilation $\dn(s)$ is strictly monotone with respect to a norm $\|\cdot\|$ in $\R^n$, 
 if there exist $\beta>0$ such that
 \begin{equation}\label{eq:strictly_monotone_dilation}
\|\dn(s)\|:=\sup_{x\neq \zero} \tfrac{\|\dn(s)x\|}{\|x\|}\leq e^{\beta s}, \quad  \forall s\leq 0.
\end{equation}
\end{Definition}
The monotonicity property implies that $\dn(s)$ acts as a strong contraction when $s<0$ and a strong expansion when $s>0$. This in turn  implies that $\forall x\in \R \backslash\{0\}$, there exists a unique pair $s_0,x_0$ such that $x = \dn(s_0)x_0$.
%The canonical homogeneous norm introduce an alternative norm topology in $\R^n$ since $\|x\|_{\dn}\leq \rho \Leftrightarrow x\in B_{\dn}(\rho)$.
  \begin{Theorem} \label{the:limit_norm_dilation}
 \cite[Corollary 6.5]{Polyakov2020:Book}  If $\dn$ is a linear dilation in $\R^n$, then the following statement holds :
 \begin{itemize}
     \item[1)]  $\dn$
 is strictly monotone with respect to the norm $\|x\|_{P}$ for $x\in\R^n$ if and only if 
\begin{equation} \label{eq:recall_monotone}
PG_{\dn}+G_{\dn}^{\top}P\succ 0, \quad P\succ 0. 
\end{equation}
\item[2)] if  $\dn$
 is strictly monotone with respect to the norm $\|x\|_{P} =\sqrt{x^\top P x}$ then 
 \begin{align}
     e^{\alpha s} \leq \|\dn(s)\|_P \leq e^{\beta s},\quad \text{if} \quad s\leq 0 \\
     e^{\beta s} \leq \|\dn(s)\|_P \leq e^{\alpha s},\quad \text{if} \quad s\geq 0
 \end{align}
 where $\alpha = \frac{1}{2}\lambda_{max}(P^{\frac{1}{2}}G_\dn P^{-\frac{1}{2}} +P^{-\frac{1}{2}}G_\dn^\top P^{\frac{1}{2}})$ and $\beta= \frac{1}{2}\lambda_{min}(P^{\frac{1}{2}}G_\dn P^{-\frac{1}{2}} +P^{-\frac{1}{2}}G_\dn^\top P^{\frac{1}{2}})$.
 \end{itemize}
\end{Theorem}

\subsection{Canonical homogeneous norm}

In the case of a strictly monotone dilation,  the so-called canonical homogeneous norm \cite{Polyakov_etal2018:TAC} can be introduced through a homogeneous projection onto the unit sphere. 
\begin{Definition}\label{def:hom_norm}\cite{Polyakov2020:Book} \itshape
Given a strictly monotone linear dilation $\dn$   with respect to the norm $\|\cdot \|_P$.
	The $\dn$-homogeneous function $\|\cdot\|_{\dn} : \R^n \mapsto [0,+\infty)$ defined as  follows
	\begin{equation}\label{eq:hom_norm}
			\|x\|_{\dn,P}:=
			\begin{cases}
				e^{s_x} : \|\dn(-s_x)x\|_P=1,\!\! &\text{if } x\neq 0,\\
				0, &\text{if } x = 0,
				\vspace{-1mm}
			\end{cases}
	\end{equation}
	is called the canonical homogeneous norm.
\end{Definition}
 As the dilation is strictly monotone then it has been demonstrated in  \cite{polyakov2019sliding} that  $\|\cdot\|_{\dn,P}$
 is single-valued, continuous everywhere on $\R^n$ and continuously differentiable on $\R^n\backslash\{0\}$
 \begin{equation}\label{eq:hom_norm_deriv}
\footnotesize
		\tfrac{\partial \|x\|_{\dn,P}}{\partial x} =\|x\|_{\dn,P} \tfrac{x^\top \dn^\top (-\ln{\|x\|_{\dn,P}})P\dn (-\ln{\|x\|_{\dn,P}}) }{x^\top \dn^\top (-\ln{\|x\|_{\dn,P}}) PG_\dn \dn (-\ln{\|x\|_{\dn,P}}) x},\,  x\neq \zero.
\end{equation} 
Moreover, for any $x\in \R^n$ and $s\in \R$, we have  the  following properties: 
 $\|x\|_{\dn,P}=\|-x\|_{\dn,P}$, $\|\dn(s)x\|_{\dn,P}=e^{s}\|x\|_{\dn,P}$ and $
\|x\|_{\dn,P}=1 \Leftrightarrow \|x\|_P=1$.
%Canonical homogeneous norm is an example of the so-called $\dn$-homogeneous function. 
\begin{Definition}
 A function $h : \R^n\to \R$ is said to be $\dn$-homogeneous of degree $\nu$ if 
 \begin{equation}
h(\dn(s)x) = e^{\nu s}h(x), \quad \forall x\in \R^n\backslash \{0\},\quad \forall s\in\R
 \end{equation}
\end{Definition}

\begin{Definition} \itshape
  \label{def:hom_operator}
     A vector field $f: \R^n \to \R^n$ is said to be $\dn$-homogeneous of degree $\nu\in \R$  if 
     \begin{equation}\label{eq:homogeneous_operator}
    	      f(\dn(s)x)=e^{\nu s}\dn(s)f(x), \quad \text{for} \quad x\in \R^n\backslash\{0\}, \forall s\in\R
         \end{equation}
\end{Definition}
The property of homogeneity of the vector field $f$ implies that  the solutions $x(t,x_0)$ of the system $\dot x=f(x), t>0, x(0)=x_0$ are also homogeneous, meaning that $x(t,\dn(s)x_0)=\dn(s) x(e^{\nu s} t, x_0)$, where $\nu\in \R$ represents the homogeneity degree of $f$.

\subsection{Homogeneous control design for linear plants}
Homogeneous control systems  offer several benefits compared to linear ones, including  faster convergence\cite{BhatBernstein2005:MCSS}, enhanced robustness \cite{Hong2001:Aut} and reduced overshoots\cite[Chapter 1]{Polyakov2020:Book}.
The following theorem recalls a procedure of the generalized homogeneous control design for linear plants and summarizes the results of the papers \cite{polyakov2016robust}, \cite{zimenko2020robust}, \cite{nekhoroshikh2021finite}.
\begin{Theorem}\label{the:homo_controller}
For the time-invariant controllable system  
\begin{equation}\label{eq:auto_sys_dis_free}
\dot{x} = A x + B u    
\end{equation}
 let a pair $\{A,B\}$ be controllable. Then
\begin{itemize}
    \item[1)] any solution $Y_0\in \R^{m\times n}, G_0\in \R^{n\times n}$ of the linear algebraic equation 
    \begin{equation}\label{eq:G0_Y0}
        AG_0 - G_0 A +BY_0 = A,\quad G_0 B =0
    \end{equation}
    is such that the matrix $G_0 - I_n$ is invertible, the matrix $G_\dn = I_n +\mu G_0$ is anti-Hurwitz for any $\mu \in [-1,\frac{1}{\tilde{n}}]$, where $\tilde{n}$ is a minimal natural number such that $\rank[B,AB,...,A^{\title{n}-1}B]=n$, the matrix $A_0 = A+ BY_0(G_0-I_n)^{-1}$ satisfies the identity 
    \begin{equation}\label{eq:A_0_homo}
        A_0 G_\dn = (G_\dn + \mu I_n) A_0,\quad G_\dn B = B
    \end{equation}
    \item[2)] the linear algebraic system 
    \begin{align}\label{eq:homo_stab}
     &\text{\footnotesize 
        $A_0 X + XA_0^\top +BY + Y^\top B + \rho (G_\dn X + XG^\top_\dn ) = 0$}\\
        &G_\dn X + XG_\dn ^\top \succ 0, \quad X = X^\top \succ 0\label{eq:G_d_mono}
    \end{align}
    has a solution $X\in \R^{n\times n}, Y\in \R^{m\times n}$ for any $\rho\in \R_+$
    \item[3)] the canonical homogeneous norm $\|\cdot\|_{\dn,P}$ induced by the weighted Euclidean norm $\|x\| = \sqrt{x^\top P x}$ with $P = X^{-1}$ is a Lyapunov function of system \eqref{eq:auto_sys_dis_free} with 
    \begin{align}
        u(x) = K_0 x + \|x\|_{\dn,P}^{\mu+1} K\dn(-\ln\|x\|_{\dn,P}) x \label{eq:pril_homo_controller}\\
        K_0 = Y_0 (G_0-I_n)^{-1},\quad K = Y X^{-1}
    \end{align}
    \item[4)] the feedback law $u$ given by \eqref{eq:pril_homo_controller} is continuously differentiable on $\R^n\backslash\{0\}$, $u$ is continuous at zero if $\mu >-1$ and $u$ is discontinuous at zeros if $\mu =-1$; 
    \item[5)] the system \eqref{eq:auto_sys_dis_free}, \eqref{eq:pril_homo_controller} is $\dn-$homogeneous of degree $\mu$.
\end{itemize}
\end{Theorem}
Clearly, when the homogeneous degree $\mu<0$, the closed-loop system \eqref{eq:auto_sys_dis_free} with controller \eqref{eq:pril_homo_controller} is uniformly finite-time stable and if $\mu>0$, it is nearly fixed-time stable. If $\mu =0$, the controller \eqref{eq:pril_homo_controller} turn out to be a linear controller $u = K_0 x + Kx$.

%{
%\color{red} The following theorem (as it is stated) is correct only for the single input case. 
%\begin{Theorem}\cite{zimenko2020robust}
%\label{the:homo_controller}
%\itshape
%For the controllable system  $\dot{x} = A x + B u$,   with nilpotent A, an  homogeneous controller can be designed in the form 
%\begin{equation} \label{eq:homo_controller}
%u(x) = \|x\|_{\dn,P}^{\mu+1} K \dn(-\ln\|x\|_{\dn,P})x
%\end{equation} 
%with $K  = YP^{-1}\in \R^{p\times n}$ and dilation $\dn$ generated by $G_\dn \in \R^{n\times n}$ satisfying 
%\begin{equation}
%AG_\dn = (G_\dn+\mu I)A,\quad
% G_\dn B =  B
%\end{equation}
%and $X = P^{-1} \in \R^{n\times n}$, $Y\in \R^{p\times n}$ under constraints 
%\begin{align}
%XA^\top + A X + Y^\top B^\top + BY &\prec 0, \\
%G_\dn P + PG_\dn^\top \succ 0,\quad X&\succ 0
%\end{align}
%\end{Theorem}
%}

\section{Homogeneous Invariant  Ellipsoid}
\label{Pre_results}
Let us consider the following MIMO (multiple-inputs multiple-outputs) non-linear system:
\begin{equation}
\label{eq:Auto_system_disturbance}
\begin{aligned}
\dot{x} &= f(x,\omega), \quad t>0
\end{aligned}
\end{equation}
where $f: \R^n \times \R^n \mapsto \R^n$ is a continuous vector field, $x$ and $\omega$ are defined as before.%, $h:\R^n\mapsto \R^p$ is a smooth mapping, and $z(t)\in \R^l$ is the controlled  output.
%,$\omega\in L^{\infty}(\R,\R^p)$ is the external bounded input (treated as a disturbance):    
%\begin{equation}\label{eq:disturbance_assumption}
% \|\omega\|_Q\leq 1,\quad Q=Q^\top \succ 0, \quad  \forall t \geq 0.
%\end{equation}

%Without loss of generality, we may assume that  $h(x)=Cx$, since there always exists a transformation of coordinate such that the later equality holds.

To have a bounded solution of \eqref{eq:Auto_system_disturbance} for any bounded perturbation $\omega$, the system must satisfy a condition like ISS\footnote{The system \eqref{eq:Auto_system_disturbance} is considered  Input-to-State Stable (ISS) \cite{Sontag1989:TAC} if there exist two functions,  $\beta$ from the class $\mathcal{KL}$ and  $\gamma$ from the class $\mathcal{K}$, such that for every bounded control input $u(\cdot)$ and initial state $x_0$, the solution satisfies $\|x(t)\|\leq \beta(\|x_0\|,t) + \gamma(\|u\|)$.}.
 %The following definition introduces the notion of an invariant ellipsoid of system  \eqref{eq:Auto_system_disturbance}.
 Recall \cite{khlebnikov2011optimization}\cite{nazin2007rejection},  that the set
 \begin{equation}
\varepsilon(P) = \{ x\in \R^n : \|x\|_P \leq 1 \}, \quad P\succ 0
\end{equation} 
is an  ellipsoid centered at origin and configured with the  matrix $P$. Similarly, we define the \textit{$\dn-$homogeneous ellipsoid} as follows
\begin{equation}
    \varepsilon_\dn(P) = \{ x\in \R^n : \|x\|_{\dn,P} \leq 1 \},\quad P\succ 0 
\end{equation}
where  $\|\cdot\|_{\dn,P}$ is the canonical homogeneous norm induced by the weighted Euclidean norm $\|\cdot\|_P$.

\begin{Definition}\label{def:ellipsoid} \itshape \cite{khlebnikov2011optimization} 
 For  system \eqref{eq:Auto_system_disturbance},\eqref{eq:disturbance_assumption},  a $\dn-$homogeneous ellipsoid $\varepsilon_{\dn}(P)$
with a configuration matrix $P$ is said to be 
\begin{itemize}
\item \textbf{Invariant}  if the condition $\forall x(0)\in \varepsilon(P)$ implies $x(t)\in \varepsilon_{\dn}(P)$ for any $t\geq 0$;
\item \textbf{Attractive} if  $x(0)\notin \varepsilon_{\dn}(P)$ then $x(t)\!\to\! \varepsilon(P)$ as $t\!\to\! \infty$.
\end{itemize}
\end{Definition}
In other words,  an  ellipsoid is considered invariant if it retains any trajectory initiated from the interior of the ellipsoid. If the ellipsoid is   attractive, any trajectory starting outside the ellipsoid will converges to (or into) it. 

The definition of the conventional invariant ellipsoid can be obtained by replacing $\varepsilon_{\dn}(P)$ with $\varepsilon(P)$ in the above definition. 
Formally, the set of the conventional ellipsoids $\varepsilon(P)$ is more rich than the set of $\dn$-homogeneous ellipsoids $\varepsilon_{\dn}(P)$, since, according to the definition of the canonical homogeneous norm, the matrix $P$ for the $\dn$-homogeneous ellipsoid must satisfy the restriction $P\in\R^{n \times n} :PG_\dn + G_\dn^\top P \succ 0$, which is required for the existence of the canonical homogeneous norm. However, for any positive definite symmetric matrix $P$,  a generator $G_{\dn}$ can always be selected such that the latter matrix inequality holds.  Indeed,  according to the  Theorem \ref{the:homo_controller} the generator can always be selected as $G_{\dn}=I_n+\mu G_0$, then for $\mu$ sufficiently close to zero we have $PG_\dn + G_\dn^\top P= 2P+\mu(PG_0+G_0^{\top} P)\succ 0$ due to the positive definiteness of $P$.
This means  by choosing an appropriate dilation $\dn$, any conventional ellipsoid can be converted into a $\dn-$homogeneous one.
Additionally, the conventional and $\dn$-homogeneous invariant ellipsoids have the same characterization.  
\begin{Lemma}\itshape\label{lemma:inv_iff}
 The following two claims are equivalent :
 \begin{enumerate}
     \item $\varepsilon_{\dn}(P)$ is a $\dn$-homogeneous  invariant ellipsoid of the system \eqref{eq:Auto_system_disturbance};
     \item $x^{\top} Pf(x,\omega)\leq 0$ for all $x \in \R^n: \|x\|_{P}=1$ and $\forall \omega\in \R^p : \|\omega\|_Q\leq 1$. 
 \end{enumerate} 
 \end{Lemma}
\textbf{Proof.}
$1)\Rightarrow 2)$  Suppose inversely that $ \varepsilon(P)$ is invariant   and  $\exists x_0 : \|x_0\|_{P}=1$, $ \exists \omega_0\in\R^p : \|\omega_0\|_Q\leq 1$ such that $x_0^{\top} Pf(x_0,\omega_0)> 0$.   Then,    for any solution $x(t)$ of the system \eqref{eq:Auto_system_disturbance} with $\omega=\omega_0$,
by the formula \eqref{eq:hom_norm_deriv} we have  
$$\tfrac{d \|x(t)\|_{\dn,P}}{dt} = r(t)\tfrac{ x^{\top}\!(t) \dn^{\top}\!(-\ln r(t)) P\dn(-\ln r(t))f(x(t),\omega_0)}{x^{\top}\!(t)\dn^{\top}\!(-\ln r(t)) PG_\dn \dn(-\ln r(t))x(t)},$$
where $r(t)=\|x(t)\|_{\dn,P}$. Since $f$ is continuous 
and $x(t)\to x_0$ as $t\to 0^+$ then 
$$\left.\tfrac{d \|x(t)\|_{\dn,P}}{dt}\right|_{t=0^+}=
\tfrac{x_0^{\top} Pf(x_0,\omega_0)}{x_0^{\top} PG_{\dn}x_0}>0,$$
where the equivalence  $\|x_0\|_{\dn,P}=1\Leftrightarrow \|x_0\|_P=1$ is utilized. The obtained inequality means that the function $t\mapsto \|x(t)\|_{\dn,P}$ is growing on some interval of time $[0,\varepsilon]$ and $\|x(t)\|_{\dn,P}>1$ for $t\in (0,\varepsilon)$. The latter contradicts to the invariance of the ellipsoid $\varepsilon(P)$. 

$2)\Rightarrow 1)$ Suppose $\varepsilon_\dn(P)$ is not invariant ellipsoid, i.e., for some $\omega$ satisfying \eqref{eq:disturbance_assumption} there exists a solution $x(t)$ of the system \eqref{eq:Auto_system_disturbance} initiated inside the ellipsoid $\varepsilon_{\dn}(P)$ such that there  exists $T>0$ : $\| x(T)\|_P >1$. Since the function $t\mapsto \|x(t)\|_{\dn}$ is continuous then  there exists $t_b \in (0,T)$ such that $\|x(t_b)\|_{\dn,P}=1$ and $1<\|x(t)\|_{\dn,P}$ for all $t\in (t_b,T)$. 
In this case, using the continuity of $f$ we derive  $\limsup_{t\to t_b} \frac{\|x(t)\|_{\dn,P}-\|x(t_b)\|_{\dn, P}}{t-t_b} =
\limsup_{t\to t_b} \frac{x^\top(t_b)P f(x(t_b),\omega(t))}{x(t_b)^{\top} PG_{\dn}x(t_b)}>0$. Hence, there exists the time instant $t^*$ (belonging to a neighborhood of $t_b$) such that $\omega(t^*)^{\top}Q\omega^*(t^*)\leq 1$ and  $\frac{x^\top(t_b)P f(x(t_b),\omega(t^*))}{x(t_b)^{\top} PG_{\dn}x(t_b)}>0$.
Therefore, we obtain the contradiction to the condition 2).  
$\blacksquare$

The invariant/attractive ellipsoid can be used to characterize the effects caused by external perturbations $\omega(t)$. For example, a conventional ellipsoid $\varepsilon(P)$ is invariant for the linear system 
\begin{equation} \label{eq:example_auto}
    \dot{x}=Ax+D\omega,\quad A\in \R^{n\times n}, D\in \R^{n\times p},
\end{equation}
where $A$ is a Hurwitz matrix, the pair $\{A,D\}$ is controllable and $x,\omega$ are defined  as before,  if and only if the LMI holds \cite{nazin2007rejection}, \cite{abedor1996linear} :
\begin{equation}
\resizebox{.43\textwidth}{!} 
{$
X^\top A^\top + AX + \beta X +\frac{1}{\beta} DQ^{-1}D^\top \preceq  0, \quad X = P^{-1}\succ 0.
$}
\end{equation}
For the considered linear system, any invariant ellipsoid is attractive \cite{nazin2007rejection}, \cite{abedor1996linear}. This statement is incorrect in general. Below we show that the same conclusion takes a place for the linear plant \eqref{eq:control_sys} with any stabilizing homogeneous controller \eqref{eq:pril_homo_controller}.

}

%\begin{Lemma}\itshape \label{lemma:invariant_attractive_necessary}
%Let $\tilde f$, $\dn$, $\dn_\omega$ be defined as in Corollary \ref{cor:nonlinear_sys_suffi_inva}.  If $\varepsilon_{\dn}(P)$ is an $\dn$-homogeneous \textbf{attractive}  ellipsoid of the system \eqref{eq:Auto_system_disturbance}, then $x^{\top} Pf(x,\omega)\leq 0$ for all $x \in \R^n: \|x\|_{\dn,P}=1$ and $\forall \omega\in \R^p : \omega^{\top}\omega\leq 1$. 
%\end{Lemma}
%\textbf{Proof.\textcolor{red}{ i think this proof is wrong...but i don't know how to correct it}}
%Suppose that $\exists x_0 : \|x_0\|_{\dn,P} =1$ and $\exists \omega_0\in \R^p : \|\omega_0\|_Q\leq 1$ such that 
%$x_0^{\top} Pf(x_0,\omega_0)> 0$ . 
%Then for system \eqref{eq:Auto_system_disturbance} with initial condition $x(0) = x_0$, we have 
%\begin{align}
%\frac{d}{dt}\|x_0\|_{\dn,P} & = \|x_0\|_{\dn,P}\frac{x_0^\top \dn^\top_x(-s)  P\dn(-s) f(x_0,\omega_0)}{x_{s0}^\top PG_{\dn} x_{s0}} \\
%& =\|x_0\|_{\dn,P}^{1+\mu} \frac{x^\top_{s0}  P f(x_{s0},\omega_{s0})}{x_{s0}^\top PG_{\dn} x_{s0}}
%> 0 \label{eq:a_i_ineq}
%\end{align}
%where $x_{s0}$ satisfying $\| x_{s0}\|_{P} =1$ with $x_{s0} = \dn(-s)x_0$ and $\|\omega_{s0}\|_Q\leq 1$ with $\omega_{s0} = \dn(-s)\omega_0$.
%\eqref{eq:a_i_ineq} means there exists an trajectory starts from the surface of ellipsoid $\|x\|_{\dn,P} =1$ will not converges to it.    This  contradicts with $\varepsilon_{\dn,P}$ being $\dn-$homogeneous attractive  ellipsoid.$\blacksquare$

\section{Homogeneous control with minimal invariant/attractive ellipsoid}
\label{Main_results}
\subsection{LMI-based characterization of invariant/attractive ellipsoid}
%The necessary and sufficient condition of a $\dn-$homogeneous ellipsoid to be invariant for the system \eqref{eq:control_sys} with the homogeneous controller \eqref{eq:pril_homo_controller} is given by the following theorem.
A criterion for determining if a $\dn-$homogeneous ellipsoid is an invariant set for the system \eqref{eq:control_sys} with the homogeneous controller \eqref{eq:pril_homo_controller} is presented in the following theorem and is both necessary and sufficient.
\begin{Theorem} \label{the:invariant_controller_desgn}\itshape
Let $D\neq 0$,  the controller $u$ be defined as in Theorem \ref{the:homo_controller} for some $X\in \R^{n\times n}$, $n\geq 3$ and some dilation $\dn$ in $\R^n$. Let $\dn_\omega$ be dilation in $\R^p$ such that the vector field 
\begin{equation}\label{eq:tilde_f}
\tilde f(x,\omega) = \begin{bmatrix}
Ax + Bu(x) + D\omega \\0
\end{bmatrix}    
\end{equation}
is $\tilde{\dn}-$homogeneous of degree $\mu$ with respect to the dilation
\begin{equation}\label{eq:dn_tilde}
    \tilde \dn(s) = \begin{bmatrix}
 \dn(s) & 0\\
 0 &\dn_{\omega}(s)\end{bmatrix}, s\in \R.
\end{equation}
in $\R^{n+p}$.
The $\dn$-homogeneous ellipsoid 
$\varepsilon_{\dn}(X^{-1})$ is \textbf{invariant} for the system \eqref{eq:control_sys}, \eqref{eq:pril_homo_controller} 
if and only if there exist $\beta>0$:  
\begin{equation}\label{eq:Invariant_iff_LMI_theo}
\resizebox{.43\textwidth}{!} 
{$
\begin{array}{c}
 W\!:=\!\left[\!\begin{array}{cc}
    A_0X\!+\!XA_0^\top  \!+\! BY\!+\!Y^\top \!B^\top \!+\!\beta  X & PD \\
    D^\top P  & -\beta Q 
    \end{array} \!\right]\! \preceq\! 0,\\
   G_{\dn}X\!+\!XG_{\dn}^{\top}\!\succ\! 0,\quad X\!\succ\! 0
\end{array}
$}
\end{equation}
with $ Y = KX \in\R^{m\times n}$.

\end{Theorem}
\textbf{Proof. } Let us denote $P=X^{-1}$.
\textit{Sufficiency}: 
From the first LMI of \eqref{eq:Invariant_iff_LMI_theo}, we derive 
\begin{equation}
\resizebox{.43\textwidth}{!} 
{$
 \begin{bmatrix}
 x \\ \omega
 \end{bmatrix}^\top  \! \begin{bmatrix}
    PA_0\!+\!A_0^\top P\!+\! PBK\!+\!K^\top B^\top P \!+\!\beta P  & PD \\
    D^\top P  & -\beta Q 
    \end{bmatrix} \begin{bmatrix}
 x \\ \omega
 \end{bmatrix}\!\preceq\!0
 $}
\end{equation}
which can be written as 
\begin{align*}
&x^\top (PA_0 + A_0^\top P+ PBK+K^\top B^\top P )x + x^\top PD \omega \\
&+ \omega^\top D^\top P x \leq \beta (\omega^\top Q \omega - x^\top P x)
\end{align*}
Since $\|x\|_{\dn,P}$ is continuous and $x\to x_0$ as $t\to 0$, 
then for $\omega^\top Q\omega\leq 1$ and $\|x_0\|_{\dn,P}=1$ we derive
\begin{align}
&\resizebox{.43\textwidth}{!} 
{$
\left.\tfrac{d \|x(t)\|_{\dn,P}}{dt}\right|_{t=0^+}=
\tfrac{x_0^{\top} (PA_0+A_0^\top P+PBK+K^\top B^\top P) x_0 + x_0^\top PD\omega + \omega ^\top D^\top P x_0}{x_0^{\top} (PG_{\dn}+G_\dn ^\top P) x_0} \leq 0
$}\\
& = \tfrac{x_0^{\top} P(A_0x_0 +Bu(x_0)+ D\omega)}{x_0^{\top} PG_{\dn} x_0} \leq 0
\end{align}
where the equivalence  $\|x_0\|_{\dn,P}=1\Leftrightarrow \|x_0\|_P=1$ and the second LMI of \eqref{eq:Invariant_iff_LMI_theo} are utilized. By Lemma \ref{lemma:inv_iff}, we conclude that 
$\varepsilon_{\dn}(P)$ is a $\dn$-homogeneous \textbf{invariant} ellipsoid of the system \eqref{eq:control_sys}.

\textit{Necessity}: 
First of all, notice that, by definition of the $\dn$-homogeneous invariant ellipsoid, the canonical homogeneous norm is well defined, which ensures that the second and the third matrix inequalities in \eqref{eq:Invariant_iff_LMI_theo}
are fulfilled.  Therefore we need to prove the fulfilement of the first one.
Let us initially show that if $\varepsilon_\dn(P)$ is a $\dn$-homogeneous \textbf{invariant} ellipsoid of the system \eqref{eq:control_sys}
then
\begin{equation}\label{eq:impl1}
 \omega^{\top} \dn^{\top}(-\ln \|x\|_{\dn,P}) Q \dn(-\ln \|x\|_{\dn,P}) \omega \leq 1 \quad \Rightarrow \quad 
\end{equation}
\begin{equation}\label{eq:impl2}
    x^\top\dn^{\top}(-\ln \|x\|_{\dn,P}) P \dn(-\ln \|x\|_{\dn,P}) f(x,\omega)\leq 0
\end{equation}

Suppose inversely that $\varepsilon_\dn(P)$ is a  $\dn-$homogeneous invariant ellipsoid, but $\exists x_1\in \R^n$, $\exists \omega_1\in \R^p$ such that 
\[
\omega^{\top}_1 \dn_\omega^{\top}(-\ln \|x_1\|_{\dn,P}) Q \dn_\omega(-\ln \|x_1\|_{\dn,P}) \omega_1 \leq 1 
\]
and
\[
    x^\top_1\dn^{\top}(-\ln \|x_1\|_{\dn,P}) P \dn(-\ln \|x_1\|_{\dn,P}) f(x_1,\omega_1)> 0.
\]
 Denoting $x_{s1} = \dn(-s) x_1$, $\omega_{s1} = \dn_\omega(-s) \omega_1$ and $s = \ln{\|x_1\|_{\dn,P}}$,  we derive
 \[
\omega^{\top}_{s1} Q \omega_{s1} \leq 1,  \quad 
     e^{\mu s} x_{s1}^\top P f(x_{s1},\omega_{s1})> 0
\]
where the homogeneity  of $\tilde{f}(x,\omega)$ is utilized on the last step.
Since $\|x_{s1}\|_P=1$ then   we derive 
 \[
 \|x_{s1}\|_P=1, \quad \omega^{\top}_{s1} Q \omega_{s1} \leq 1,  \quad 
     x_{s1}^\top P f(x_{s1},\omega_{s1})>0.
\]
The latter contradicts Lemma \ref{lemma:inv_iff}, so the implication \eqref{eq:impl1}-\eqref{eq:impl2} takes a place.

Moreover, using the notation $\pi= \dn(-\ln \|x\|_{\dn,P}) x$ and $\omega_x=\dn_{\omega}(-\ln \|x\|_{\dn,P}) \omega$,  
we derive 
$\dn(-\ln \|x\|_{\dn,P}) f(x,\omega)=e^{\mu s}f(\pi, \omega_x)=e^{\mu s}[(A_0+BK)\pi+D\omega_x]$.
Since for any $x\neq 0$ we have $\pi^{\top} P \pi=1$ then 
the proven implication can be rewritten 
as follow :
\[
 \omega^{\top}_{x} Q \omega_x \leq 1\quad \text{ and }\quad \pi^{\top}P\pi=1 \quad \Rightarrow \quad 
\]
\[
    \pi^\top P [(A_0+BK)\pi+D\omega_x]\leq 0.
\]

%Obviously there exists a pair $\bar \pi $ and $\bar\omega_x$ such that $\bar\omega^{\top}_{x} Q \bar\omega_x < \bar\pi^{\top} P\bar\pi$.
%This together with the claim above allows to
Applying Lemma \ref{lemma:S_procedure_new} we derive 
$\exists \tau_1\geq 0, \exists \tau_2\in\R$: $\tau_1+\tau_2\leq 0$, 
\[
\left[
\begin{smallmatrix}
P(A+BK)+(A+BK)^{\top}P & PD\\
D^{\top}P & 0
\end{smallmatrix}
\right]\leq \tau_1\left[
\begin{smallmatrix}
0 & 0\\
0 & Q
\end{smallmatrix}
\right]+\tau_2 \left[
\begin{smallmatrix}
P & 0\\
0 & 0
\end{smallmatrix}
\right].
\]
If the latter inequality holds for $\tau_2<-\tau_1$ then it holds for $\tau_2=-\tau_1$.
Indeed, $\tau_2\leq -\tau_1\; \Rightarrow \; \tau_{2}P\preceq -\tau_1 P$. 
Taking into account $P=X^{-1}$ we conclude that 
the inequality \eqref{eq:Invariant_iff_LMI_theo} holds for $\beta:=\tau_1\geq 0$. 
Let us show that $\beta>0$.  Indeed, if $\beta=0$ then the following inequality must hold : 
\begin{equation}\label{eq:zero_schur}
 \left[
     \begin{smallmatrix}
    A_0X\!+\!XA_0^\top  \!+\! BY\!+\!Y^\top B^\top & D \\
    D^\top  & 0 
    \end{smallmatrix}\right] \! \preceq\! 0. 
\end{equation}
The latter inequality may hold only if $D=\zero$. Indeed,  by the generalized Schur complement \cite{zhang2006schur}, it has the equivalent conditions  $A_0X\!+\!XA_0^\top  \!+\! BY\!+\!Y^\top B^\top\preceq 0$ and $I\cdot D^\top=0$.  
Therefore, $\beta$ may be zero only if $D=0$. This contradicts the assumption of the theorem.
$\blacksquare$

Similarly to the case of the linear control system \cite{nazin2007rejection}, \cite{abedor1996linear}, the invariant ellipsoids of the $\dn$-homogeneous control system \eqref{eq:control_sys}, \eqref{eq:pril_homo_controller} can be characterized by means of a linear matrix inequality \eqref{eq:Invariant_iff_LMI_theo}, as long as the perturbations $\omega$ are involved into the system in a generalized homogeneous manner. The only difference in the LMI \eqref{eq:Invariant_iff_LMI_theo} with respect to the conventional case is the condition $G_{\dn}X\!+\!XG_{\dn}^{\top}\!\succ\! 0$, which, in the view of Theorem \ref{the:homo_controller}, disappears as $\mu\to 0$. So, the condition of the homogeneity of the vector field $\tilde f$ is the main difference.
The mentioned condition introduces a restriction to the structure of the exogenous perturbations, namely, to  a class of admissible  matrices $D$ in the system \eqref{eq:control_sys}.

In the general case, the homogeneity condition can be checked as follows. 
\begin{Proposition}\label{pro:G0_omega}
Let $G_0$  and the homogeneous   controller $u$ be defined as in Theorem  \ref{the:homo_controller}.
If there exists a matrix $G_{0,\omega}\in \R^{p\times p}$ such that
\begin{equation}\label{eq:G0_omega}
    G_0D=DG_{0,\omega}
\end{equation}
then for any $\mu \in \R$ satisfying 
\begin{equation}\label{eq:mu_G0_omega_1}
1+\mu +\mu\mathrm{Re}\left(\lambda_i(G_{0,\omega})\right)> 0, \quad \forall i=1,2,...,n
\end{equation}
\begin{enumerate}
\item the matrix $G_{\dn_\omega}=I_p+\mu (I_p+G_{0,\omega})$ is anti-Hurwitz and $\dn_{\omega}(s)=e^{s G_{\dn_\omega}},s\in \R$, is a dilation in $\R^p$;
\item the vector field $\tilde f$ defined by \eqref{eq:tilde_f} is $\tilde \dn$-homogeneous of degree $\mu$, where $\tilde \dn$ is a dilation in $\R^{n+p}$ defined by \eqref{eq:dn_tilde}.
\end{enumerate}
\end{Proposition}

\textbf{Proof. }  The inequality  \eqref{eq:mu_G0_omega_1} immidiately implies that $G_{\dn_\omega}$ is anti-Hurwitz and  $\dn_\omega(s) =e^{s G_{\dn_\omega}} $  is a dilation in $\R^p$.
From equation \eqref{eq:G0_omega}, we derive that 
\begin{equation}
     G_0^{i}D=DG_{0,\omega}^{i}
\end{equation}
where $i$ is an integer and $i\geq 0$.
Then we have
\begin{equation}
\resizebox{.43\textwidth}{!} 
{$
 e^{\mu G_0 s}D=\sum_{i=0}^{+\infty} \frac{G_{0}^i \mu^i s^i}{i!}D=  D\sum_{i=0}^{+\infty} \frac{G_{0,\omega}^i \mu^i s^i}{i!}=De^{\mu G_{0,\omega} s}
$}
\end{equation}
The latter implies
\begin{equation}
  e^{\mu s} \dn(s) D = D\dn_{\omega}(s) 
\end{equation}
Finally, since $G_0$ is defined as in Theorem \ref{the:homo_controller}, then we have  
\[
    f(\dn(s)x,\dn_\omega(s)\omega) =A\dn(s)x+Bu(\dn(s)x)+D\dn_{\omega}(s)\omega=
    \]
    \[
    e^{\mu s} \dn(s)x+e^{\mu s} \dn(s)Bu(x)+e^{\mu s} \dn(s)D\omega=e^{\mu s} \dn(s) f(x,\omega),
\]
where $f(x,\omega) = Ax + Bu(x) + D\omega$.
%which implies $\tilde{f}(\dn(s) x,\dn_\omega(s) \omega) = e^{\mu s } \tilde{\dn}(s) \tilde{f}(x,\omega)$.
$\blacksquare$

It is well-know \cite{nazin2007rejection}, \cite{abedor1996linear} that, in many cases, an invariant ellipsoid of the linear control system  is attractive as well. The same conclusion can be made for the considered homogeneous control system in the case when the dilation $\dn_{\omega}$ in $\R^p$ is strictly monotone with respect to the norm $\|\cdot\|_Q$ in $\R^p$.  
\begin{Corollary} \label{cor:iff_attractive}
Under the conditions of Theorem \ref{the:invariant_controller_desgn}, the invariant ellipsoid $\varepsilon_{\dn}(X^{-1})$ is attractive provided that 
\begin{equation}\label{eq:monotone_omega}
QG_{\dn_{\omega}}+G_{\dn_{\omega}}^{\top} Q\succ 0.
\end{equation}
\end{Corollary}

\textbf{Proof. } Considering the canonical homogeneous norm as a Lyapunov function of the closed-loop system  \eqref{eq:control_sys}, \eqref{eq:pril_homo_controller},  we derive
\[
\frac{d}{dt}\|x\|_{\dn,P}=\|x(t)\|_{\dn,P}^{1+\mu } \tfrac{x_s^\top(t) X^{-1} [(A_0+BK)x_s(t)+D\omega_s]}{x_s^\top(t) PG_{\dn} x_s(t)}
\]
\[
=\tfrac{\|x(t)\|_{\dn,P}^{1+\mu }}{x_s^\top(t) PG_{\dn} x_s(t)}\left(\left[\!\begin{smallmatrix}
 X^{\text{--}1}x_s \\ \omega_s
 \end{smallmatrix}\!\right]^{\!\top}  \! W \left[\begin{smallmatrix}
 X^{\text{--1}}x_s \\ \omega_s
 \end{smallmatrix}\right]\!-\beta +\beta \omega_s^{\top}\!Q \omega_s\right),
\]
where $x_s=\dn(-\ln \|x\|_{\dn,P})x$ and $\omega_s=\dn(-\ln \|x\|_{\dn,P})\omega$. 
Since \[
\omega^{\top}\dn_{\omega}^{\top}(s)Q\dn_{\omega}(s)\omega+
\int^0_s \frac{d}{d\tau}\omega^{\top}\dn_{\omega}^{\top}(\tau)Q\dn_{\omega}(\tau)\omega  d\tau =
\omega^{\top}Q\omega  
\]
then, due to \eqref{eq:monotone_omega}, for any $s>0$ and any $\omega\neq 0$ we have 
\[
\int^0_s \frac{d}{d\tau}\omega^{\top}\dn_{\omega}^{\top}(\tau)Q\dn_{\omega}(\tau)\omega d\tau=
\]
\[
=\int^0_s \omega^{\top}\dn_{\omega}^{\top}(\tau)\left(QG_{\dn_{\omega}}+G_{\dn_{\omega}}^{\top} Q\right)\dn_{\omega}(\tau)\omega d\tau>0, 
\]
and $
\omega^{\top}\dn_{\omega}^{\top}(s)Q\dn_{\omega}(s)\omega< \omega^{\top}Q\omega$.
 The latter implies that $\omega^{\top}_s Q \omega_s<\omega^{\top}Q\omega\leq 1$ for 
any $x:\|x\|_{\dn,P}>1$ and any $\omega\neq 0$. Therefore, using the inequality \eqref{eq:Invariant_iff_LMI_theo}  with $\beta>0$, 
we conclude that 
\[
\frac{d}{dt}\|x\|_{\dn,P}<0, \quad\forall x:  \|x\|_{\dn,P}>1, \forall \omega : \|\omega\|_{Q}\leq 1. 
\]
The latter implies that the $\dn$-homogeneous ellipsoid is attractive. $\blacksquare$
\begin{Corollary}\label{cor:mu_G0_omega}
If $G_{0,\omega}\in \R^{p\times p}$ is defined as in Proposition \ref{pro:G0_omega},
then the inequalities \eqref{eq:mu_G0_omega_1} and \eqref{eq:monotone_omega} with  $G_{\dn_{\omega}}=I_p+\mu(I_p +G_{0,\omega})$ hold for any $\mu \in \R$ satisfying 
\begin{equation}\label{eq:mu_G0_omega}
\resizebox{.43\textwidth}{!} 
{$
1+\mu+\mu\tfrac{\lambda_{i}\left(Q^{\frac{1}{2}} G_{0,\omega}Q^{-\frac{1}{2}}\!+\!Q^{-\frac{1}{2}} G_{0,\omega}^{\top}Q^{\frac{1}{2}}\right)}{2}\!>\!0, \, \forall i=1,...,n.
$}
\end{equation}
\end{Corollary}
\textbf{Proof.} 
Since 
\[
1+\mu+\mu\tfrac{\lambda_{i}\left(Q^{\frac{1}{2}} G_{0,\omega}Q^{-\frac{1}{2}}\!+\!Q^{-\frac{1}{2}} G_{0,\omega}^{\top}Q^{\frac{1}{2}}\right)}{2}=
\]
\[
\tfrac{\lambda_{i}\left(2(1+\mu)I_n+\mu Q^{\frac{1}{2}} G_{0,\omega}Q^{-\frac{1}{2}}\!+\!Q^{-\frac{1}{2}} G_{0,\omega}^{\top}Q^{\frac{1}{2}}\right)}{2}=
\]
\[
\tfrac{\lambda_{i}\left(Q^{\frac{1}{2}} G_{\dn_\omega}Q^{-\frac{1}{2}}\!+\!Q^{-\frac{1}{2}} G_{\dn_\omega}^{\top}Q^{\frac{1}{2}}\right)}{2}
\]
then the inequality \eqref{eq:mu_G0_omega} implies the inequality \eqref{eq:monotone_omega}. Taking into account that $Q\succ 0$, by Lyapunov inequality, we 
derive that $G_{\dn_\omega}$ is anti-Hurwitz, i.e., the inequality \eqref{eq:mu_G0_omega_1} holds as well.
$\blacksquare$

Obviously, the inequalities \eqref{eq:mu_G0_omega_1} and \eqref{eq:mu_G0_omega} are always feasible, at least, for $\mu$ close to $0$.
So, the existence of  solution for the equation \eqref{eq:G0_omega} is the only critical restriction of an applicability of the homogeneous invariant/attractive ellipsoids method.

\subsection{Minimal invariant and attractive ellipsoids}
Inspired by \cite{nazin2007rejection} and Theorem \eqref{the:invariant_controller_desgn}, the  optimal tuning of the homogeneous controller \eqref{eq:pril_homo_controller} can be formulated in terms of the following Semi-Definite-Programming (SDP)  problem: 
 \begin{equation}\label{eq:obj_fun}
     \mathrm{tr}(X) \quad \rightarrow \quad min
 \end{equation}
subject to the LMI constraints \eqref{eq:Invariant_iff_LMI_theo}.

%Since minimize the  $\dn-$homogeneous ellipsoid $\varepsilon_{\dn,P} : \|x\|_{\dn,P}\leq 1$ is equivalent to minimize the ellipsoid $\varepsilon(P) : \|x\|_P\leq 1$ and the output of system is proportional to the system state. Here we select the trace of matrix to characterize the size of output ellipsoid which requires to minimize $tr[X]$. 

Note that for any fixed $\beta$, the problem of finding an optimal solution reduces to  minimizing the linear function \eqref{eq:obj_fun} subject to the LMI constraints    \eqref{eq:Invariant_iff_LMI_theo}, which is the classical SDP problem.  There
exist many MATLAB toolboxes for its numerical solution such as SeDuMi and YALMIP.
The considered SDP problem \textit{without} the additional constraint $G_{\dn}X+XG_\dn^{\top}\succ 0$ was studied in \cite{nazin2007rejection}, where  its  convexity has been proven. %Since the optimal solution $X$ of this convex optimization problem  only  exists for  $\beta$ belonging to a bounded interval.  Therefore the optimal $X$ can be found as well after adding one more constrains of $X$. 

Notice that the linear controller is a particular case of the homogeneous controller. Indeed, by Theorem \ref{the:homo_controller},
 for $\mu=0$ we have $\dn(s)=e^s I_n$, $\dn(s)=e^{s}I_{p}$,
\begin{equation}
    u = K_0x + K x. \label{eq:linear_controller}
\end{equation}
and  the condition of the homogeneity of the vector field $\tilde f$ is always fulfilled (see Theorem \ref{the:invariant_controller_desgn}).
In this case, the LMI \eqref{eq:Invariant_iff_LMI_theo} becomes 
\begin{equation}
\resizebox{.43\textwidth}{!} 
{$
 A_0X+XA_0^\top +\beta X + BY+Y^\top B^\top + \tfrac{1}{\beta} DQ^{-1}D^\top \preceq 0, X\succ 0 \label{eq:Linear_op_LMI_S} 
 $}
\end{equation}
which perfectly fits the results   of \cite{nazin2007rejection}, \cite{abedor1996linear}. Moreover, if the system of algebraic equations \eqref{eq:G0_Y0}, \eqref{eq:G0_omega} has a solution then an optimal (in the sense of minimal invariant ellipsoid) linear 
controller can be upgraded to a (nonlinear) homogeneous one without any degradation of the minimal invariant/attractive ellipsoid.  
\begin{Corollary} 
 Let the tuple $(X_{opt},Y_{opt},\beta_{opt})\in \R^{n\times n} \times \R^{m\times n} \times (0,+\infty)$ be a solution of the SDP problem 
\eqref{eq:obj_fun}, \eqref{eq:Linear_op_LMI_S} then for any $\mu\in \R$ satisfying
\begin{equation}\label{eq:mu_G0}
1+\mu\lambda_{\min}\left(X^{-\frac{1}{2}}_{opt} G_{0}X^{\frac{1}{2}}_{opt}+X^{\frac{1}{2}} G_{0}^{\top}X^{-\frac{1}{2}}_{opt}\right)>0,
\end{equation}
the  tuple $(X_{opt},Y_{opt},\beta_{opt})$ is a solution of the SDP problem 
\eqref{eq:obj_fun}, \eqref{eq:Invariant_iff_LMI_theo},
with $G_{\dn}=I_n+\mu G_{0}\in \R^{n \times n}$. Moreover,  
the $\dn$-homogeneous ellipsoid $\varepsilon_{\dn}(X^{-1})_{opt}$ is 
\begin{itemize}
\item \textbf{invariant} for  the homogeneous control system \eqref{eq:control_sys}, \eqref{eq:pril_homo_controller} provided that $\mu\in \R$ satisfies \eqref{eq:mu_G0} and \eqref{eq:mu_G0_omega_1};
\item \textbf{attractive} for the homogeneous control system \eqref{eq:control_sys}, \eqref{eq:pril_homo_controller} provided that $\mu\in \R$ satisfies \eqref{eq:mu_G0} and \eqref{eq:mu_G0_omega}. 
\end{itemize}
\end{Corollary}
\textbf{Proof. } 
The required result directly  follows from the equivalence of inequality   \eqref{eq:mu_G0}  to $G_\dn X + X G_\dn ^\top \succ 0$.
$\blacksquare$

A transformation of a linear optimal controller to a homogeneous optimal controller can be useful, for example,  if an additional criterion needs to be optimized or an additional condition has to be satisfied. For example,  if $K_0=0$ and the selection $\mu=-1$ is admissible in the latter corollary then the corresponding ("upgraded") homogeneous controller is uniformly bounded. Indeed, 
\[
u^{\top}u= x^{\top}\dn^{\top}(-\ln\|x\|_{\dn,P})K^{\top}K\dn(-\ln\|x\|_{\dn,P}) x\leq 
\]
\[
\lambda_{\max}(P^{-1/2}K^{\top}KP^{-1/2})
\]
where the identity $x^{\top}\dn^{\top}(-\ln\|x\|_{\dn,P})P\dn(-\ln\|x\|_{\dn,P}) x=1$ is utilized in  the last step.
\begin{Corollary} 
If $K_0=0$ and $\mu=-1$ then the homogeneous control \eqref{eq:pril_homo_controller} satisfies the restriction
\[
u^{\top}(x)u(x)\leq \bar u^2,\quad \quad  \forall x\in \R^n,
\]
for some $\bar u>0$, if and only if 
\begin{equation}\label{eq:LMI_baru}
\left[
\begin{array}{cc}
X &Y^{\top}\\
Y & \bar u^2 I_{m}
\end{array}\right]\succeq 0
\end{equation}
where $X = P^{-1} ,Y = KX$.
\end{Corollary}
\textbf{Proof. } 
If $K_0 = 0$ and $\mu =-1$, we have $$u^\top u \leq  \lambda_{max}(P^{-1/2}K^{\top}KP^{-1/2})$$
The latter inequality is sharp, thus the inequality $u^{\top}u\leq \bar{u}^2 $ is equivalent to $\lambda_{max}(P^{-1/2}K^{\top}KP^{-1/2})\leq \bar{u}^2$ and to
\begin{align}
    X^\top K^\top KX \leq \bar{u}^2 X,\quad  
\end{align}
Apply the  Schur complement we derive \eqref{eq:LMI_baru}.
$\blacksquare$

Therefore, the design of globally bounded controllers that minimize the size of the invariant or attractive ellipsoid can be achieved by using the $\dn$-homogeneous invariant ellipsoid method. The design procedure in this case is formulated in terms of the 
SDP problem \eqref{eq:obj_fun}, \eqref{eq:Invariant_iff_LMI_theo}, \eqref{eq:LMI_baru}.

\section{Simulation and Experiment}\label{sim_results}

\subsection{Numerical Results}
In this subsection, the linearized model of rotary inverted pendulum system Quanser Qube Servo-2 is studied.
\begin{align}
\dot{x} &= A x + Bu +D \omega \label{eq:ex_state}
%y &= C x +E u  \label{eq:ex_output}
\end{align}
where 
%\begin{align*}
%A  &= \begin{psmallmatrix}
%0& 1& 0 &0\\
%0 & 0 & 1 & 0\\
%1& 0 & 0 & 1\\
%0 & 0 & 0 &1
% \end{psmallmatrix},\quad 
% B = \begin{psmallmatrix}0 \\0\\0\\1 \end{psmallmatrix},\quad
%  D = \begin{psmallmatrix}
% 0 &0\\ 
% 0&1\\
% 0&0\\
% 0.1& 1 \end{psmallmatrix} 
% C &= \begin{psmallmatrix}
% 0& 0&0&0 \\
% 0& 1& 0 & 0
% \end{psmallmatrix},\quad
% E = \begin{psmallmatrix}1\\0\end{psmallmatrix}
%,\quad x_0 = \begin{bmatrix}
%1.5\\ -1 \\ 0 \\  0
%\end{bmatrix}
%\end{align*}
%Note that the system above is controllable and $E^\top C = 0$.
\begin{align*}
%\resizebox{.495\textwidth}{!} 
%{$
A  &= \begin{psmallmatrix}
0& 0& 1 &0\\
0 & 0 & 0 & 1\\
0& \frac{l^2rg m_p^2}{J_t} & \frac{-b_rJ_p}{J_t}-\frac{k_m^2}{R_m}\frac{J_p}{J_t} & \frac{-l r m_pb_p}{J_t}\\
0 & \frac{g l m_p  J_r}{J_t}  & -\frac{l r m_p  b_r}{J_t}-\frac{k_m^2}{R_m}\frac{l r m_p }{J_t} & \frac{-J_r b_p}{J_t}
\end{psmallmatrix},\\
B &=\frac{k_m}{R_m} \begin{psmallmatrix}
0\\ 0\\ \frac{J_p}{J_t}\\ \frac{l r m_p }{J_t}
\end{psmallmatrix},
%$}
\quad
 D = \begin{psmallmatrix} 
      1 & 0& 0& 0 \\ 
      0 &1 & 0 &0  \\
      0 &0 & 1 &0\\
      0 &0 &0 &1
      \end{psmallmatrix},\quad x = \begin{psmallmatrix}
          x_1 \\ x_2 \\x_3 \\x_4
      \end{psmallmatrix} =\begin{psmallmatrix}
          \theta\\ \alpha\\ \dot{\theta}\\ \dot{\alpha}
      \end{psmallmatrix},
\end{align*}
$\theta$ is the angle position of arm and $\alpha$ is the angle position of pendulum.
\begin{table}
 \caption{Pendulum Model Parameters }
 \label{pendulum_parameter_table}
 \begin{center}
 \begin{tabular}{|c|c|c|c|}
 \hline
  Parameter & Description & Value & Units\\
 \hline\hline
 $R_m$ & Motor Resistance & 8.4 & $\Omega$ \\
$K_m$  & Back-emf constant & 0.042 & $\frac{V\cdot s}{rad}$ \\
 $m_r$ & Rotary arm mass & 0.095 & kg\\ 
 $r$ & Rotary arm length & 0.085 & m\\ 
 $J_r$ &Rotary Inertia Moment & $\frac{r^2 m_r}{3}$ & $kg\cdot m^2$ \\
 $b_r$ &  Rotary Damping Coefficient & $10^{-3} $& $ \frac{N\cdot M\cdot s}{rad}$\\
 $m_p$ & Pendulum Link Mass  & $0.024$ & $kg$ \\
 $L_p$ & Pendulum Link length & $0.129$ & $m$  \\
 $l$ & Pendulum center of mass & $\frac{l_p}{2}$ & $m$\\
 $J_p$ &Pendulum inertia moment & $ \frac{m_p\cdot L_p^2}{3}$ & $kg\cdot m^2$ \\
 $b_p$ & Pendulum Damping Coefficient & $5\times 10^{-5}$& $ \frac{N\cdot M\cdot s}{rad}$  \\
 $g$ & Gravity constant & 9.81 & $\frac{m}{s^2}$\\
 \hline
 \end{tabular}
 \end{center}
 \end{table}
The external disturbance is bounded by $\|\omega\|_Q\leq 1$ with 
\begin{align*}
Q = \begin{psmallmatrix} 2& 0 & 0 & 0\\
                    0 & 2 & 0 &0 \\
                    0 & 0 & 1 & 0\\
                    0 & 0 & 0 & 2
                    \end{psmallmatrix},\quad \omega = \begin{bmatrix} 0.2&0.3&0.3&0.4\end{bmatrix}^\top \cdot \sin{\frac{t}{2}}
\end{align*}
The first step is using Theorem \ref{the:homo_controller} and Proposition \ref{pro:G0_omega}  to find the matrix $Y_0$, $G_0$ and $G_{0,\omega}$ such that \eqref{eq:G0_Y0} and \eqref{eq:G0_omega} hold.
\begin{align}
    Y_0 &= \begin{bmatrix}
  0  & 10.65 &   -0.73 &   0.47
    \end{bmatrix},\\
G_0 & = \begin{psmallmatrix}
    -3 & 2.02 & 0 & 0\\
    0 & -1 & 0 & 0 \\
    0 & 0.38 & -2 & 2.02\\
    0 & 0 & 0 & 0
\end{psmallmatrix}, G_{0,\omega} =
\begin{psmallmatrix}
    -3 & 2.02 & 0 & 0\\
    0 & -1 & 0 & 0 \\
    0 & 0.38 & -2 & 2.02\\
    0 & 0 & 0 & 0
\end{psmallmatrix}
\end{align}
Then $K_0$ and dilation $G_\dn$ is found by $Y_0$ and $G_0$.

Then the second step is to find the optimal linear controller \eqref{eq:linear_controller} with \eqref{eq:obj_fun}  \eqref{eq:Linear_op_LMI_S}. 
\begin{equation}
    K = \begin{bmatrix}
    27.12& -177.13 & 10.91 & -17.93
    \end{bmatrix}
\end{equation}
In the third step, the optimal linear gain K is applied to defined $Y=KX$ in the optimal problem \eqref{eq:obj_fun} under constraint \eqref{eq:homo_stab} and \eqref{eq:G_d_mono}.
\begin{equation}
    X = \begin{psmallmatrix}
    1.33 & 0.11 & -0.87 & 0.42 \\
    0.11 & 0.05 & -0.51 & -0.58 \\
    -0.87 & -0.51 & 48.52 &35.47 \\
    0.42 & -0.58 & 35.47 & 30.13
    \end{psmallmatrix}
\end{equation}
This optimal $X$ obtained above is needed to define the homogeneous controller \eqref{eq:pril_homo_controller}.

\begin{figure}[ht]    
\centering
\includegraphics[scale=0.2]{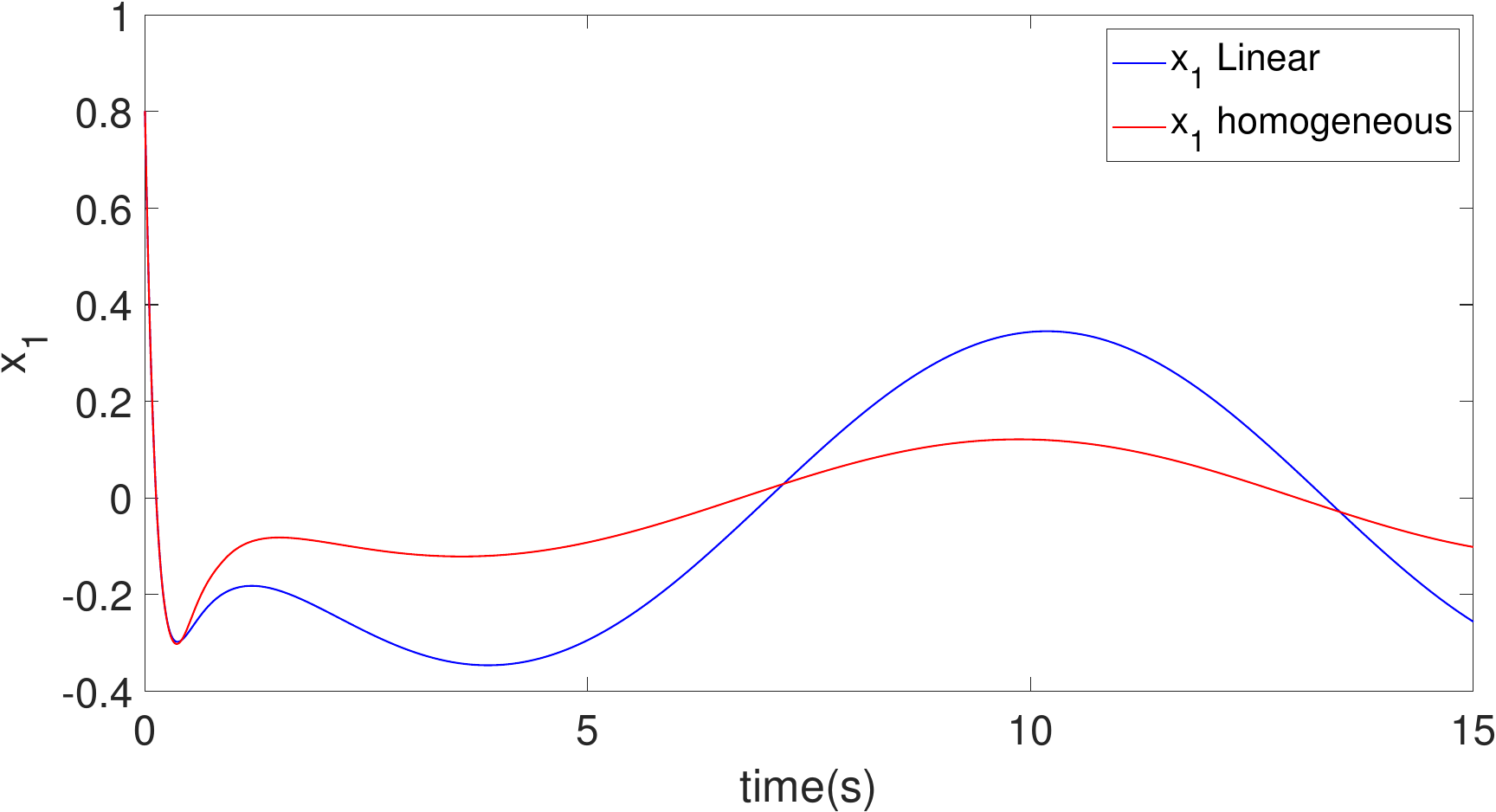}  
\caption{Comparison of  $x_1$ by  linear and homogeneous controller}
\label{fig:x_1_compare07}
\end{figure}

\begin{figure}[ht]    
\centering
\includegraphics[scale=0.2]{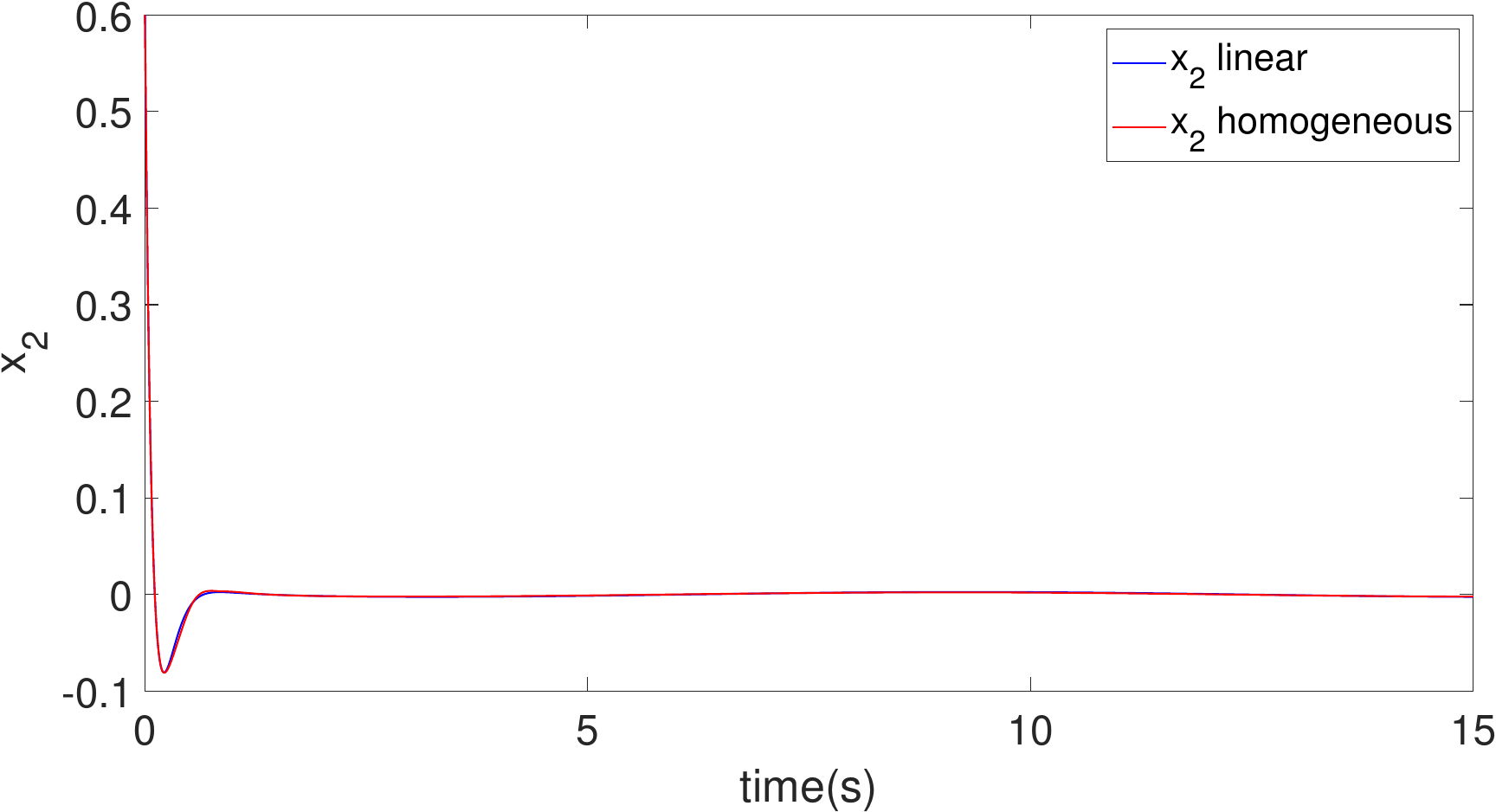}  
\caption{ Comparison of  $x_2$ by  linear and homogeneous controller}
\label{fig:x_2_compare07}
\end{figure}

\begin{figure}[ht]    
\centering
\includegraphics[scale=0.2]{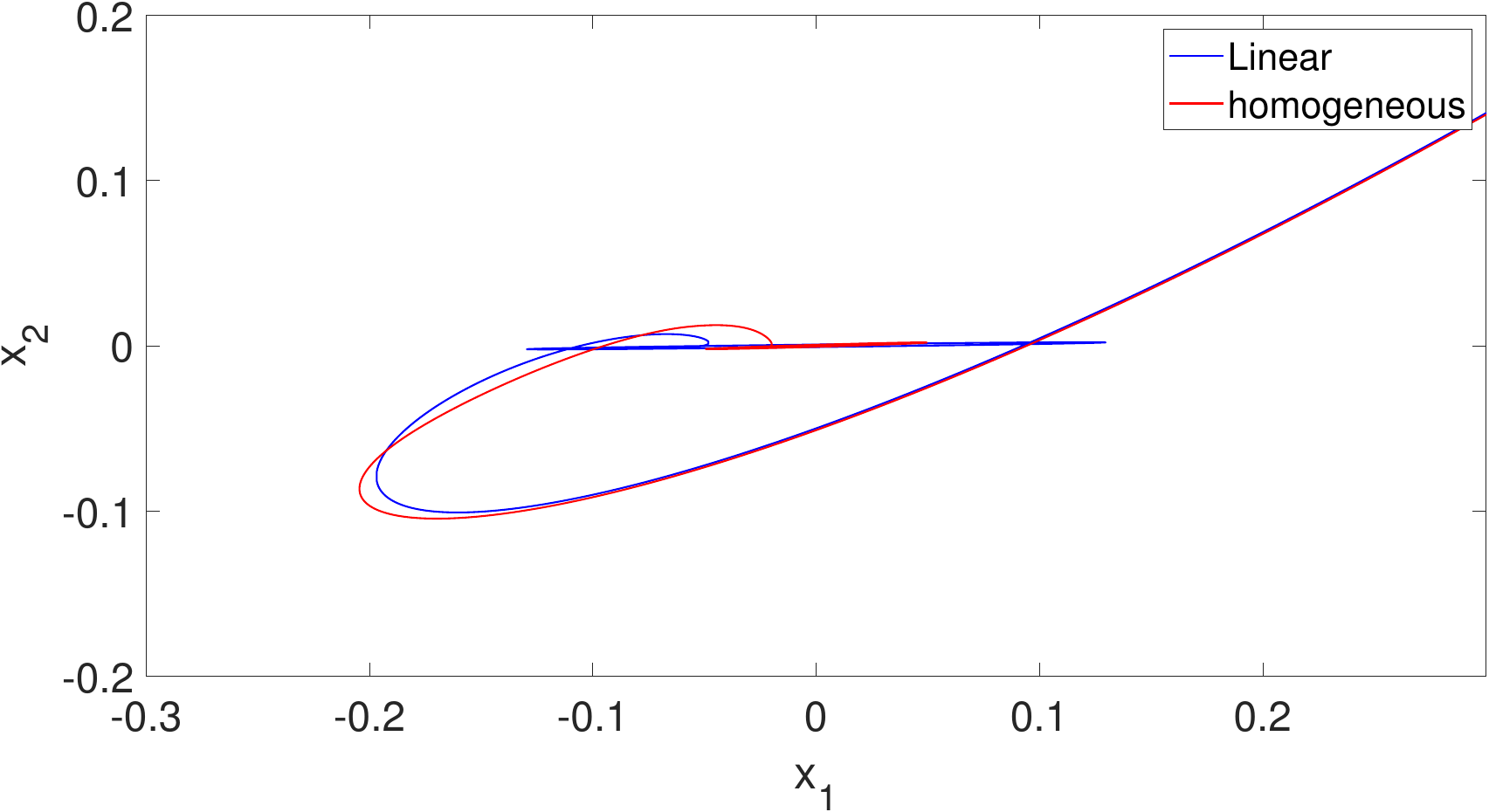}  
\caption{Comparison of trajectory $[x_1,x_2]$ by linear and homogeneous controller  }
\label{fig:Trajectory07}
\end{figure}

In the final step, we can verify that for $\mu =-0.7$, both groups  \eqref{eq:control_sys} \eqref{eq:pril_homo_controller} and \eqref{eq:mu_G0}  \eqref{eq:mu_G0_omega} are feasible.
The simulation results are  shown in Fig.\ref{fig:x_1_compare07}-\ref{fig:Trajectory07}. According to the Fig. \ref{fig:x_1_compare07}, for linear controller the maximum value of  $x_1$ in  the invariant set is about $0.1296$. The same $x_1$ by using homogeneous controller is about $0.0494$. The precision improvement on $x_1$ is about $60\%$. Besides, we can also find  that $x_1$ has a faster response by using the homogeneous controller.

%\begin{figure}[ht]    
%\centering
%\includegraphics[scale=0.2]{Pictures/x_3_Comparison07.eps}  
%\caption{ Comparison of  $x_3$ by  linear and homogeneous controller}
%\label{fig:x_3_compare07}
%\end{figure}

%\begin{figure}[ht]    
%\centering
%\includegraphics[scale=0.2]{Pictures/x_4_Comparison07.eps}  
%\caption{ Comparison of  $x_4$ by  linear and homogeneous controller}
%\label{fig:x_4_compare07}
%\end{figure}

\subsection{Experiment results}
This experiment is based on the rotary inverted pendulum  Quanser Qube Servo-2, where the same model \eqref{eq:ex_state} is applied. The system parameters  are presented in Table \ref{pendulum_parameter_table}. Here we consider more physical model of disturbances $D = \begin{bmatrix}0& 0 &2.53 &2.50
\end{bmatrix}^\top $ and $Q=1$.  We use the well tuned linear feedback gain $K = \begin{bmatrix}2 &-35 &1.5& -3 \end{bmatrix}^\top$ provided by Quanser to design the homogeneous controller by upgrading the linear one.   We obtain 
\begin{equation}
    X = \begin{psmallmatrix}
   0.00159 &  0.00002 & -0.00275 & -0.00018\\
   0.00002 &   0.00002 & 0.00011& -0.00004\\
  -0.00275 &  0.00011  & 0.02492&   0.01377 \\
  -0.00018 & -0.00004  & 0.01377 &  0.01210
    \end{psmallmatrix}
\end{equation}
for $\mu =-0.7$ satisfying  \eqref{eq:control_sys} \eqref{eq:pril_homo_controller} and \eqref{eq:mu_G0}  \eqref{eq:mu_G0_omega}.
 
The experiment is repeated five times to guarantee a fair comparison. The results  are presented in the Table \ref{tab:linear}  and Table \ref{tab:homo}. It is clear to see that the stabilization precision of $\theta$ and $\alpha$ in $L^2$-norm is improved by using the  homogeneous controller about $20.22\%$ and $8.3\%$,  respectively. The stabilization precision in the $L^{\infty}$-norm is improved  about $11.62\%$ and $15.4\%$, respectively.
The energy consumption increases about $2.41\%$.

\begin{table}[h!]
\centering
\caption{Pendulum stabilization with linear controller }
\label{tab:linear}
\begin{tabular}{||c |c|c| c|c|c ||}
 \hline
 Linear &$\|\theta\|_{L^\infty}$& $\|\alpha\|_{L^\infty}$& $\|\theta \|_{L^2}$&$\|\alpha \|_{L^2}$ & $\|u\|_{L^2}$ \\ [0.5ex]
 \hline\hline
Test 1 & 0.0767&0.0092 &0.21333 & 0.01347 & 0.79136  \\ 
 \hline
Test 2 &	0.07363& 0.00613 & 0.20377  &	0.01271 &	0.7717   \\
 \hline
Test 3 &	0.07363&0.00613&0.21414 &	0.01282 &	0.67597  \\
 \hline
Test 4 &	0.0859&0.0092 & 0.2379 &	0.01423 &	0.74302  \\
 \hline
Test 5 &	0.0859&0.0092 & 0.24353 &	0.01332 &	0.60481 \\ 
 \hline\hline
 Average  & 0.0792& 0.00797 &    0.2225 &   0.0133 & 0.7174\\
 \hline
\end{tabular}
\end{table}

%\begin{table}[h!]
%\centering
%\caption{Pendulum stabilization with homogeneous controller }
%\label{tab:homo}
%\begin{tabular}{||c |c| c| c||} 
% \hline
% Homogeneous & $\|e_\theta \|_{L^2}$ & $\|e_\alpha\|_{L^2}$ & $\|u\|_{L^2}$ \\ [0.5ex]
% \hline\hline
%Test 1 & 0.15032 & 0.01422& 0.99976   \\ 
% \hline
%Test 2 & 0.15216 &0.01374& 0.99736  \\
% \hline
%Test 3 & 0.15106 &0.01321& 1.051   \\
% \hline
%Test 4 & 0.14497 &0.01334& 1.09393   \\
% \hline
%Test 5 & 0.1437 &0.01322& 1.05151  \\ 
% \hline\hline
% Average & 0.1484 & 0.013566  & 1.0387 \\
% \hline
%\end{tabular}
%\end{table}
%The optimal X applied in the homogeneous controller is found by $\mu=-0.7$ and %$D =I\in\R^{4\times 4}$.

\begin{table}[h!]
\centering
\caption{Pendulum stabilization with homogeneous controller }
\label{tab:homo}
\begin{tabular}{||c |c| c|c| c|c||} 
 \hline
 Homogeneous & $\|\theta\|_{L^\infty}$& $\|\alpha\|_{L^\infty}$& $\|\theta \|_{L^2}$ & $\|\alpha\|_{L^2}$ & $\|u\|_{L^2}$\\ [0.5ex]
 \hline\hline
Test 1 &	0.0675&0.00613& 0.18762 & 0.01268 & 0.6292 \\ 
 \hline
Test 2&	0.0675 &0.00613& 0.18695 &	0.01277 &	0.77291   \\
 \hline
Test 3&	0.0675&0.00613 & 0.16906 &	0.01126 &	0.71843    \\
 \hline
Test 4 &	0.0767& 0.0092& 0.17481 &	0.01245 &	0.76564    \\
 \hline
Test 5 &	0.07056 &0.00613& 0.16912 &	0.01202 &	0.78741  \\ 
 \hline\hline
 Average & 0.0700& 0.00674 &0.1775 & 0.0122  & 0.7347 \\
 \hline
\end{tabular}
\end{table}
 
%Fig. \ref{fig:theta_compare}-\ref{fig:input_compare} show a comparison plot of Test 3 with linear and homogeneous controller.
%\begin{figure}[ht]    
%\centering
%\includegraphics[scale=0.16]{Pictures/theta_Comparison_zoom.eps}  
%\caption{Trajectory of  $\theta$ by  linear and homogeneous controller}
%\label{fig:theta_compare}
%\end{figure}
%
%\begin{figure}[ht]    
%\centering
%\includegraphics[scale=0.25]{Pictures/alpha_Comparison.eps}  
%\caption{Trajectory of  $\alpha$ by  linear and homogeneous controller}
%\label{fig:alpha_compare}
%\end{figure}
%
%\begin{figure}[ht]    
%\centering
%\includegraphics[scale=0.25]{Pictures/Input_Comparison.eps}  
%\caption{Trajectory of linear and homogeneous controller  $u(t)$ }
%\label{fig:input_compare}
%\end{figure}
 
%%%%%%%%%%%%%%%%%%%%%%%%%%%%%%%%%%%%%%%%%%%%%%%%%%%%%%%%%%%%%%%%%%%%%%%%%%%%%%%%%

\section{Conclusion}\label{sec:conclusion}
This article extends the invariant/attractive ellipsoid method \cite{nazin2007rejection,khlebnikov2011optimization,Poznyak_etal2014:Book} to a class of  the generalized homogeneous system. The LMI-based characterization of $\dn$-homogeneous invariant/attractive ellipsoid for a linear plant is obtained by the homogeneous control. It shows that, under certain restriction to the structure of perturbations, an optimal (in the sense of the minimal invariant ellipsoid) $\dn-$homogeneous invariant/attractive ellipsoid can be upgraded from  any conventional invariant/attractive one  by properly selecting the dilation $\dn$. Despite that theoretically optimal invariant ellipsoids are the same for both linear and homogeneous controllers, the simulation and experiment results show that the optimal homogeneous controller may provide a better precision one under the same conditions/perturbations.

%This article extends the invariant ellipsoid method to the design of homogeneous optimal controller.  Since in the real platform an optimized linear controller is easier to obtain than the nonlinear one, a procedure of upgrading a linear controller  to a generalized homogeneous one while preserving the minimal invariant  set is proposed. The simulation results show that the nonlinear bounded  homogeneous controller is able to  further minimize the effect of disturbance. 

%This article firstly proposed some equivalent conditions of $\dn$-homogeneous ball being an invariant set and then presents the approach of  designing an homogeneous  controller based on the invariant set method. The proposed  homogeneous  controller based on the optimized linear one can guarantee the same minimal invariant set to the linear case.

%This article proposed a simple approach to design a globally bounded homogeneous  controller that has the same minimal invariant set to linear case.  The optimized   linear controller is obtained  by the method of invariant ellipsoid, since the concept of invariant set allows one to reduce the design of optimal controller to the problem in terms of LMIs and convex minimization.  The proposed  globally bounded   homogeneous  controller is proven feasible to  further minimize the effect of disturbance with the optimal linear controller, which is supported by the simulation results of two-mass-spring system. 

% References
%\section*{Acknowledgment}

%This work is supported by the Project Medibot.

\appendix
\begin{Lemma}[Proposition 4.1, \cite{Polyak1998:JOTA}]\label{lemma:S_procedure_new}
Let $n\geq 3$. Let us consider the homogeneous quadratic forms $f_i(x) = x^\top A_i x$, $i=0,1,2$ with $x\in \R^n$ $A_i = A_i^\top\in \R^{n\times n}$, and the numbers $\alpha_0,\alpha_1,\alpha_2\in \R$, $\alpha_2\neq 0$. Let there exist $\mu_1,\mu_2\in \R$ and $x_0\in \R^n$ such that 
$$\mu_1 A_1+\mu_2 A_2\succ 0 \text{ and }f_1(x_0)<\alpha_1, f_2(x_0)=\alpha_2.$$ 
Then the following two claims are equivalent:
\begin{itemize}
    \item[1)]$\exists \tau_1\!\geq 0, \exists \tau_2\!\in \R$ : $
    A_0 \!\preceq\! \tau_1A_1+\tau_2 A_2, \;\; \alpha_0 \!\geq\! \tau_1\alpha_1+\tau_2 \alpha_2$;
\item[2)] 
$f_0(x)\!\leq\! \alpha_0, \;\; \forall x\in \R^n:    f_1(x) \!\leq \! \alpha_1,\; f_2(x) \!=\! \alpha_2 $.
\end{itemize}
\end{Lemma}

\bibliographystyle{IEEEtran}
\bibliography{Invariant_set_Ref}\ %IEEEabrv instead of IEEEfull

\end{document}